\documentclass[12pt]{amsart}
\usepackage[utf8]{inputenc}
\usepackage{amsmath}
\usepackage{xcolor}
\usepackage{amssymb}
\newtheorem{theorem}{Theorem}
\newtheorem*{theorem*}{Theorem}

\newtheorem{lemma}{Lemma}

\newtheorem*{acknowledgements*}{Acknowledgements}

\makeatletter
\def\blfootnote{\gdef\@thefnmark{}\@footnotetext}
\makeatother

\def\house#1{\setbox1=\hbox{$\,#1\,$}%
\dimen1=\ht1 \advance\dimen1 by 2pt \dimen2=\dp1 \advance\dimen2 by 2pt
\setbox1=\hbox{\vrule height\dimen1 depth\dimen2\box1\vrule}%
\setbox1=\vbox{\hrule\box1}%
\advance\dimen1 by .4pt \ht1=\dimen1
\advance\dimen2 by .4pt \dp1=\dimen2 \box1\relax}

\begin{document}
\title{A short note  on the divisibility of class numbers of real quadratic fields}
\author{Jaitra Chattopadhyay}
\address{Harish-Chandra Research Institute, HBNI, Chhatnag Road, Jhunsi, Allahabad - 211019, INDIA}
\email[Jaitra Chattopadhyay]{jaitrachattopadhyay@hri.res.in}
\begin{abstract}
For any integer $l\geq 1$, let   $p_1, p_2, \ldots, p_{l+2}$ be distinct prime numbers $\geq 5.$ For all real numbers $X>1,$ we let $N_{3,l}(X)$ denote the number of   real quadratic fields $K$ whose absolute discriminant $d_K\leq X$ and $d_K$ is divisible by $(p_1\ldots p_{l+2})$ together with  the class number $h_K$ of $K$ divisible by $2^{l}\cdot 3.$ Then, in this short note, by  following the method in  \cite{Byeonkoh},  we prove that  $N_{3,l}(X) \gg X^\frac{7}{8}$ for all large enough $X$'s.
\end{abstract}
\maketitle

\section{Introduction}

The problem of the divisibility of class numbers of number fields has been of immense interest to number theorists for quite a long time. Many mathematicians have studied the divisibility problem of class number for quadratic fields.

\smallskip

Nagell \cite{nagell} showed that there exist infinitely many imaginary quadratic fields whose class numbers are divisible by a given positive integer $g.$ Later, Ankeny and Chowla \cite{AC} also proved the same result. The  analogous result for real quadartic fields had been proved by Weinberger \cite{berger}, Yamamoto \cite{moto} and many others. 

\smallskip

Apart from the qualitative results, a great deal of work has also been done towards the quantitative versions. For a positive integer $g,$  we let  $N^+_g (X)$ (respectively, $N^-_g (X)$) denote the number of real (respectively, imaginary) quadratic fields $K$ whose absolute value of the discriminant is $d_K\leq X$ and the class number $h_K$ is divisible by $g.$ Then the general problem is to find lower bounds for the magnitude of $N^+_g (X)$ (respectively, $N^-_g (X)$) as $X \rightarrow \infty.$ M. Ram Murty \cite{Ram} proved that, for any integer $g \geq 3,$ the inequalities $N^-_g (X) \gg X^{\frac{1}{2}+\frac{1}{g}}$ and $N^+_g (X) \gg X^\frac{1}{2g}$ hold. The behavior of $N^+_g (X)$ and $N^-_g (X)$ have also been studied in \cite{Luca}, \cite{Sound} and \cite{Yu}.  

\smallskip

The case $g=3$ has been studied in  \cite{Byeonkoh}, \cite{Kalyan}, \cite{Dav}, \cite{Honda}, \cite{KishiMiyake} and in many other papers. Byeon \cite{Byeon5} studied the case $g=5$ and $g=7$ and in both the cases, he showed that $N^+_g (X) \gg X^\frac{1}{2},$ which is an improvement over the main result of \cite{Ram}.

\smallskip

In this short note, we shall study the following related problem. For a given natural number $l\geq 1,$ we fix $l+2$ distinct prime numbers $\geq 5,$ say, $p_1, \ldots, p_{l+2}$ and let $g\geq 3,$ $g\ne p_i$ for all $i,$ be a given odd integer. We let 
$N_{g,l}(X)$ denote the number of real quadratic fields $K$ whose absolute discriminant $d_K\leq X,$ $h_K$ is divisible by $g$ and $d_K$ is divisible by $(p_1\ldots p_{l+2}).$ This, in turn, implies that $2^l$ divides $h_K.$ Then the general problem is to find  the magnitude of $N_{g, l}(X).$  

\smallskip

In this short note, by adopting the method of \cite{Byeonkoh}, we prove the following theorem.
 
\begin{theorem}\label{thm1}
We have, $N_{3, l}(X) \gg X^\frac{7}{8}$ for all large enough real numbers $X.$ 
\end{theorem}

\section{Preliminaries}

In \cite{KishiMiyake}, Kishi and Miyake gave a complete classification of quadratic fields $K$ whose class number $h_K$ is divisible by $3$ 
as follows:

\smallskip

\begin{lemma} \label{lem1}
Let $g(T)=T^3 - uwT - u^2 \in \mathbb{Z}[T]$ be a polynomial with integer coefficients $u$ and $w$  such that $gcd(u,w)=1,$ $d=4uw^3 - 27u^2$ is not a perfect square in $\mathbb{Z}$ and one of the following conditions holds:
\begin{itemize}
\item[(i)]$3 \nmid w$
\item[(ii)]$3 \mid w, uw \not\equiv  3 \pmod 9, \text { and } u \equiv w \pm 1 \pmod 9 $
\item[(iii)]$3 \mid w, uw \equiv 3 \pmod 9, \text { and } u \equiv w \pm 1 \pmod {27} $
\end{itemize}
If $g(T)$ is irreducible over $\mathbb{Q},$ then the roots of  the polynomial $g(T)$ generate an unramified cyclic cubic extension $L$ over $K:=\mathbb{Q}(\sqrt{d})$ (which in turn implies, by Class Field Theory, that $3$ divides $h_K).$ Conversely, suppose $K$  is a quadratic field over $\mathbb{Q}$ with $3$ dividing the class number $h_K.$ If $L$ is an  unramified cyclic cubic extension over $K,$ then $L$  is  obtained by adjoining the roots of $g(T)$ in $K$ for some suitable choices of $u$ and $w.$
\end{lemma}

\bigskip

Using Lemma \ref{lem1},   
Byeon and Koh  \cite{Byeonkoh} proved the following result.

\begin{lemma}\label{lemma2}
Let $m$ and $n$ be two relatively prime positive integers satisfying $m \equiv 1 \pmod {18}$  and  $n \equiv 1 \pmod {54}.$
  If the polynomial $f(T)=T^3 - 3mT - 2n$ is irreducible over $\mathbb{Q},$ then the class number of the quadratic field $\mathbb{Q}(\sqrt{3(m^3 - n^2)})$ is divisible by $3.$
\end{lemma}

In \cite{Sound}, Soundararajan proved the following result (see also \cite{Byeonkoh}).

\begin{lemma}\label{lem3}
Let $X$ be a large positive real number and $T=X^\frac{1}{16}.$ Also, let $M= \frac{T^\frac{2}{3}X^\frac{1}{3}}{2}$ and  $N=\frac{TX^\frac{1}{2}}{2^4}.$ If  $N(X)$ denotes  the number of positive square-free integers $d \leq X$ with at least one integer solution $(m, n, t)$ to the equation 
\begin{equation}\label{eq1}
m^3 - n^2=27t^2d
\end{equation}
satisfying  $T < t \leq 2T, M < m \leq 2M, N < n \leq 2N, gcd(m,t)= gcd(m,n)= gcd(t,6)=1, m\equiv 19 \pmod {18 \cdot 6} \text{ and } n \equiv 55 \pmod {54 \cdot 6},$  then, we have,
$$
N(X) \asymp \frac{MN}{T} + o(MT^{\frac{2}{3}}X^{\frac{1}{3}}) \gg X^{7/8}.
$$
\end{lemma}

The following result was proved  in \cite{Kalyan} which provides a lower bound of the number of irreducible cubic polynomials with bounded coefficients.

\begin{lemma}\label{lem4}
Let $M$ and $N$ be two positive real numbers. Let 
\begin{eqnarray*}
\mathcal{S}&=&\left\{f(T) =T^3+mT+n\in\mathbb{Z}[T] \  : \  |m| \leq M, |n| \leq N, f(T) \mbox{ is irreducible over } \mathbb{Q}\right. \\
&& \left.\qquad\qquad\qquad\qquad\qquad \qquad \mbox{ and } D(f) = -(4m^3+27n^2) \mbox{ is not a perfect square}\right\}
\end{eqnarray*} 
be a subset of $\mathbb{Z}[T].$
Then $|\mathcal{S}| \gg MN.$
\end{lemma}

\section{Proof of Theorem \ref{thm1}}

Let  $l\geq 1$ be an integer and let $g=2^l \cdot 3.$  Let $K$ be  a real quadratic field over  $\mathbb{Q}$ and its class number is $h_K.$ Note that  $h_K \equiv 0 \pmod g$ if and only if $h_K \equiv 0 \pmod {2^l} $ and $h_K \equiv 0 \pmod 3.$

\bigskip

\noindent{\bf Claim 1.}  The number of quadratic field $\mathbb{Q}(\sqrt{d})$ satisfying $d\leq X,$  $d$ is divisible by $p_1p_2\ldots p_{l+2}$ and $3$ divides $h_K$ is $\gg X^\frac{7}{8}.$ 

\bigskip

For each $i = 1, 2, \ldots, l+2,$  let $a_i$ and $b_i$ be integers such that 
\begin{equation} \label{eq2}
3a_i - 2b_i \not\equiv 0 \pmod {p_i}.
\end{equation}
Then,  consider the simultaneous congruences
\begin{eqnarray*}
X&\equiv& 19\pmod{18 \cdot 6}\\
X  &\equiv& 1 + a_i p_i \pmod{p_i^2},
\end{eqnarray*}
for all $i = 1, 2, \ldots, l+2.$ Then,  by the Chinese Reminder Theorem, there is a unique integer solution $m$  modulo $18 \cdot 6 \displaystyle\prod_{i=1}^{l+2}p_i^2.$ Thus, the number of such integers $m\leq X$ is $\left((1+o(1)\right) X/\left(18 \cdot 6 \displaystyle\prod_{i=1}^{l+2}p_i^2\right)$ as $X \rightarrow \infty.$ Let $N_1(X)$ be the set of all such integers $m\leq X.$   Similarly, we consider the simultaneous congruences 
\begin{eqnarray*}
X&\equiv& 55\pmod{54 \cdot 6}\\
X  &\equiv& 1 + b_i p_i \pmod{p_i^2},
\end{eqnarray*}
for all $i = 1, 2, \ldots, l+2.$  Then,  by the Chinese Reminder Theorem, there is a unique integer solution $n$  modulo $54 \cdot 6 \displaystyle\prod_{i=1}^{l+2}p_i^2.$  Thus, the number of such integers $n\leq X$ is $\left((1+o(1)\right) X/\left(54 \cdot 6 \displaystyle\prod_{i=1}^{l+2}p_i^2\right)$ as $X \rightarrow \infty.$ Let $N_2(X)$ be the set of all such integers $n\leq X.$ 

\bigskip

Let $X$ be a large positive real number and $T=X^\frac{1}{16}.$ Also, let 
$$
\displaystyle M= \frac{T^\frac{2}{3}X^\frac{1}{3}}{2} \mbox{ and } \displaystyle N=\frac{TX^\frac{1}{2}}{2^4}.
$$ 
Now, we shall count the number of tuples $(m, n, t)$ satisfying \eqref{eq1} 
with  $T < t \leq 2T,$  $M < m \leq 2M,$ $N < n \leq 2N,$ gcd$(m,t)=$ gcd$(m,n)=$ gcd$(t,6)=1,$ $m \in N_1(X),$  $n \in N_2(X)$ with square-free integer $d.$ Then, by Lemma \ref{lem3}, we see that 
\begin{equation}\label{eq4}
N(X) \gg X^{7/8}.
\end{equation}
Now note that for any integers $m\in N_1(X)$ and $n\in N_2(X)$, we see that 
$$
m^3 - n^2 = (a_ip_i+1)^3 - (b_ip_i+1)^2 \equiv 3a_ip_i + 1 - 2b_ip_i -1 \equiv p_i(3a_i -2b_i) \pmod{p_i^2},
$$
for all $i = 1, 2, \ldots, l+2.$ By \eqref{eq2}, since $3a_i -2b_i \not\equiv 0\pmod{p_i}$ for all $i,$ we see that $m^3 - n^2 \not\equiv 0 \pmod{p_i^2}$ and hence $p_i$ divides the square-free part of $m^3-n^2$ which is $d.$ Thus, $p_1p_2\ldots p_{l+2}$ divides $d$ for all such  $d$'s. 

\bigskip

In order finish the proof of Claim 1, by Lemma \ref{lemma2}, it is enough to count the number of tuples $(m,n,t)$ satisfying \eqref{eq1}  for which $f(X) = T^3 - 3mT -2n$ is irreducible over $\mathbb{Q}.$ By Lemma \ref{lem4}, the number of such irreducible polynomial $f(T)$ is at least $\gg MN \gg X^\frac{7}{8},$ which proves Claim 1.

\bigskip

Now, to finish the proof of the theorem, we see that at least $\gg X^{7/8}$ number of real quadratic fields $K=\mathbb{Q}(\sqrt{d})$ satisfying $h_K\equiv 0\pmod{3}$ and $\omega(d) \geq l+2,$ where $\omega(n)$ denotes the number of distinct prime factors of $n.$ Therefore, by Gauss' theory of genera, we conclude that the class number $h_K$ of corresponding real quadratic field is divisible by $2^{l}$ also. Combining this fact with Claim 1, we get the theorem.

\begin{acknowledgements*} 
I would like to sincerely thank Prof. Florian Luca for his insightful remarks and having fruitful discussions with him. I acknowledge Prof. R. Thangadurai for going through the article and suggesting valuable points to improve the presentation of the paper. I am thankful to the referee for his/her valuable suggestions to improve the presentation of the paper. I am grateful to Dept. of Atomic Energy, Govt. of India and Harish-Chandra Research Institute for providing financial support to carry out this research.
\end{acknowledgements*}


\begin{thebibliography}{9999}


\bibitem{AC} 
N. Ankeny and S. Chowla,  On the divisibility of the class numbers of quadratic fields, {\it Pacific J. Math.}, {\bf 5} (1955),  321-324.

\bibitem{Byeon5}
D. Byeon,  Real quadratic fields with class number divisible by $5$ or $7,$ {\it Manu. Math.},  {\bf 120} (2006)  (2) 211-215.

\bibitem{Byeonkoh}
D. Byeon and E. Koh, Real quadratic fields with class number divisible by $3,$ {\it Manu.  Math.}, {\bf 111}  (2003), 261-263.

\bibitem{Kalyan}
K. Chakraborty and M. Ram Murty, On the number of real quadratic fields with class number divisible by $3,$ {\it Proc.  Amer. Math. Soc.}, {\bf 131} (2002), 41-44.

\bibitem{Dav}
H. Davenort and H. Heilbronn, On the density of discriminants of cubic fields, {\it Proc. Royal Soc. A}, {\bf 322}  (1971), 405-420.

\bibitem{Honda}
T. Honda, On real quadratic fields whose class numbers are multiples of $3,$ {\it J.Reine Angew. Math}, {\bf 223} (1968), 101-102.



\bibitem{KishiMiyake}
Y. Kishi and K. Miyake, Parametrization of the quadratic fields whose class numbers are divisible by three, {\it J. Number Theory}, {\bf 80} (2000), 209-217.

\bibitem{Luca}
F. Luca, A note on the divisibility of class numbers of real quadratic fields, {\it C. R. Math. Acad. Sci. Soc. R. Can}, {\bf 25} (2003) (3), 71-75.

\bibitem{nagell}  
T. Nagell, Uber die Klassenzahl imaginar quadratischer Zahkorper, {\it Abh. Math. Seminar Univ. Hamburg}, {\bf 1}  (1922) (1), 140-150.

\bibitem{Ram}
M. Ram Murty, Exponents of class groups of quadratic fields, {\it Topics in Number theory (University Park, PA, (1997) Math. Appl.},  {\bf 467} {\it Kluwer Acad. Publ., Dordrecht},  (1999), 229-239.



\bibitem{Sound}
K. Soundararajan, Divisibility of class numbers of imaginary quadratic fields, {\it J. London Math. Soc.}, {\bf 61} (2000) (2), 681-690.

\bibitem{berger} 
P. Weinberger, Real quadratic fields with class numbers divisible by $n,$  {\it J. Number Theory},  {\bf 5}  (1973),  237-241.

\bibitem{moto}
Y. Yamamoto, On unramified Galois extensions of quadratic number fields, {\it Osaka J. Math.}, {\bf 7} (1970), 57-76.

\bibitem{Yu}
G. Yu, A note on the divisibility of class numbers of real quadratic fields, {\it J. Number Theory}, {\bf 97} (2002), 35-44.


\end{thebibliography}
\end{document}